\documentclass{raex}
\usepackage{amsmath, amssymb, amsthm, amsfonts}
\usepackage[dotinlabels]{titletoc}
\usepackage{version} 

\titlecontents{section}[1.0em] 
{\bfseries}
{\contentslabel{1.0em}}
{\hspace*{-2.3em}}
{\titlerule*[0.7pc]{.}\contentspage}
\usepackage{buczofnxsurvey}
\numberwithin{figure}{section}
\begin{Author} 
	\FirstName{Zolt\'an}\LastName{Buczolich}
	\PostalAddress{Department of Analysis, ELTE E\"otv\"os Lor\'and
University, P\'azm\'any P\'eter S\'et\'any 1/c, 1117 Budapest, Hungary}
	\Email{zoltan.buczolich@ttk.elte.hu }
	\Thanks{The project leading to this application has received funding from the European Research Council (ERC) under the European Union's Horizon 2020 research and innovation programme (grant agreement No. 741420).
This author was also supported by the Hungarian National Research, Development and Innovation Office--NKFIH, Grant 124003 and  at the time of part of the work on this paper was holding a visiting researcher position 
at the R\'enyi Institute.}
\end{Author}

\begin{MathReviews} 
	\primary{28A20} 
	\secondary{40A05}
	\secondary{60F15}
	\secondary{60F20}
\end{MathReviews}

\begin{KeyWords} 
\end{KeyWords}

\title{Almost everywhere convergence questions of series of translates of non-negative functions}

  \newtheorem{theorem}{Theorem}[section]
 \newtheorem{corollary}[theorem]{Corollary}
 
 \newtheorem{lemma}[theorem]{Lemma}
 \newtheorem{proposition}[theorem]{Proposition}
 \newtheorem{example}[theorem]{Example}
 {\theoremstyle{definition}
 \newtheorem{definition}[theorem]{Definition}

 \newtheorem{question}{Problem}
 \newtheorem{problem}[theorem]{Question}
 \newtheorem{conjecture}[theorem]{Conjecture}
 }

 \newcommand{\N}{\ensuremath{\mathbb N}} 
\newcommand{\T}{\ensuremath{\mathbb T}} 
\newcommand{\R}{\ensuremath{\mathbb R}} 
\newcommand{\Q}{\ensuremath{\mathbb Q}} 
\newcommand{\Z}{\ensuremath{\mathbb Z}} 
\newcommand{\beq}{\begin{equation}}
\newcommand{\eeq}{\end{equation}}

\begin{document}

\maketitle
\begin{center}
{\sc Dedicated to the memory of Jean-Pierre Kahane}
\end{center}

\begin{abstract}
This survey paper is based on a talk given at the 44th Summer Symposium in Real Analysis in Paris.

This line of research was initiated by a question of Haight and
Weizs\"aker concerning almost everywhere convergence properties of series of the form
$\sum_{n=1}^{{\infty}}f(nx)$.

A more general, additive version of this problem is the following:

 Suppose $\Lambda$ is a discrete infinite set of nonnegative real numbers.
We say that $ {\Lambda}$ is of type 1 if the series $s(x)=\sum_{\lambda\in\Lambda}f(x+\lambda)$
 satisfies a zero-one law. This means that for any non-negative measurable 
 $f: {\ensuremath {\mathbb R}}\to [0,+ {\infty})$ either the convergence set $C(f, {\Lambda})=\{x: s(x)<+ {\infty} \}= {\ensuremath {\mathbb R}}$ modulo sets of Lebesgue zero, or its complement the divergence set $D(f, {\Lambda})=\{x: s(x)=+ {\infty} \}= {\ensuremath {\mathbb R}}$ modulo sets of measure zero.
 If $ {\Lambda}$ is not of type 1 we say that $ {\Lambda}$ is of type 2.
 
 The exact characterization of type $1$ and type $2$ sets is still not known.
 
 The part of the paper discussing results concerning this question is based on several joint papers written at the beginning with J-P. Kahane and D. Mauldin, later with B. Hanson, B. Maga and G. V\'ertesy. 

 Apart from results from the above project we also cover historic background, other related results and open questions.
 
\end{abstract}

\setcounter{tocdepth}{2}
\tableofcontents

\section{Introduction}\label{*secintro}
The 44th Summer Symposium in Real Analysis took place exactly 
5 years after the passing away of Jean-Pierre Kahane. It was partly held at
 Université Paris-Saclay, formerly named Université Paris-Sud.
Jean-Pierre Kahane was the second president of this university from 1975 to 1978.
My talk at the symposium and this survey paper are about a project in which he also participated.

This line of research started for me in the Summer of 1998 when there was a conference held in Miskolc, Hungary. During the trainride to Miskolc I discussed
some Mathematics with Dan Mauldin and he suggested that we look at an
unsolved problem from 1970, originating
from the Diplomarbeit of
 Heinrich von Weizs\"aker \cite{[HW]}.

\begin{problem}\label{*prWei}
{\it Suppose $f:(0,+ {\infty})\to [0,+\infty)$ is a measurable function.
Is it true that
$\sum_{n=1}^{{\infty}}f(nx)$
either converges (Lebesgue) almost everywhere or diverges almost everywhere, i.e.
is there  a zero-one law for $\sum f(nx)$?}
\end{problem}

In  \cite{[BM1]}
we answered the Haight--Weizs\"aker problem.
Jean-Pierre Kahane was also interested in this question, so soon afterwards 
we started to work
together with him on related questions \cite{[BKM1]}, \cite{[BKM2]} and \cite{[BKM2err]}. Recently with Bruce Hanson and two graduate students, Bal\'azs Maga and G\'asp\'ar V\'ertesy we continued this line of research 
in papers
\cite{[BMV]},
\cite{[BHMVr]},
and \cite{[BHMVa]}.

This paper also contains some open problems. If a problem is open it is listed as a problem and is numbered increasingly throughout the whole paper, 
if it is a problem/question for which the answer is known it is listed as a question and is numbered sectionwise. (In Section 
\ref{*secfnx} there are three 
problems quoted from the American Mathematical Monthly with their original numbering,
these are of course not open problems.)
Some relevant open problems from the literature are also listed as problems. 
For them we tried to do a literature check about their status and are listed only  if we were unable to find any evidence that they have been solved. 

First, in Section \ref{*secfnxp} we discuss some results related to $\sum_{n=1}^{{\infty}}f(nx)$ when $f$ is periodic, especially the ones about 
Khinchin's conjecture.

In Section \ref{*secfnx} the non-periodic case is discussed
along with some background history of questions leading to 
Question \ref{*prWei}. 
At the end of this section our solution to 
this question is stated in Theorem \ref{*HWth} along with a reminder 
of Kronecker's Theorem.

Section \ref{*sect12} contains the additive version and generalization of the
Haight--Weizs\"aker problem which lead to the definition of type $1$ and $2$
sets.
One can think of type $1$ sets as "simple/good" sets, since  for the 
corresponding series there is a zero-one law. Type $2$ sets are the "complicated/bad" sets for which there are witness functions which 
demonstrate that there is no zero-one law.
This section also contains Theorem \ref{*BKMth5}, 
a theorem,  from a joint paper, \cite{[BKM1]} with J-P. Kahane and D. Mauldin.  
This theorem gives a sufficient condition for type $2$. It was later used many times.

Section \ref{*sectdyad} discusses  constructions of some type $1$ and $2$ sets,
using dyadic rationals. We also mention some variants of these constructions
and the loosely connected question about the density of $\{ (3/2)^n \}$.

In Section \ref{*sectwit} we give two answers to questions asked in \cite{[BKM1]}.
These questions were about some structural properties of type $2$ sets and witness functions.
The first question concerned the witness function, whether it can always be selected
to be a characteristic function.
The other question is about the possible structure of convergence/divergence sets for type $2$ sets.

Looking at Examples \ref{*exdyada}, \ref{*exdyadb} and Theorem \ref{mkthm} one might think that adding more 
elements to a "good" type $1$ set can lead to a "bad" type $2$ set.
In Section \ref{*secsubsup} one can see that the hierarchy of containment of type $1$ and $2$ sets is much more complicated. For example according to Theorem \ref{sandwich} one can find a nested seqence of sets for which every other is type $1$ 
(and alternate ones are type $2$).

In the final Section \ref{*seccnxpln} we look at the case of weighted series.

\section{$\sum f(nx)$, the periodic case}\label{*secfnxp}

In this section we discuss results in the case where the function $f$ is periodic.
Without limiting generality we suppose that it has period $1$.

If instead of the series $\sum f(nx)$
we consider only the sequence
 $f(nx)$
then Mazur and Orlicz \cite{[MO]}
proved in 1940
that
$$\limsup_{n\to\oo}f(nx)=\hbox{ess}\sup f$$
for almost every
$x$. This property holds in the more general case if the
sequence
$(n)$ is replaced by an arbitrary sequence
$(\lll_{n})$ converging to $\oo.$
\medskip

This implies that if the periodic function
$f$ is the characteristic function of a set of
 positive measure,
then for almost every $x$ we have
$\sum_{n}f(nx)=\oo$.
Thus in the periodic case we have a zero-one law.

Considering
Ces{{$\grave{\hbox{a}}$}}ro $1$ means of the partial sums
of our series we are led to

\begin{conjecture}[Khinchin \cite{[Kh]} (1923)]

{\it Assume that
$E {\subset} (0,1)$
is a measurable set and
$f(x)=\chi_{E}(\{x\})$,
where $\{x\}$ denotes the fractional part of $x$.
Is it true that for almost every $x$
\begin{equation}\label{*khi}
\frac{1}{k}\sum_{n=1}^{k} f(nx)\to  {\mu}(E)?
\end{equation}}
(In our paper $ {\mu}$ denotes the Lebesgue measure.)

\end{conjecture}

This conjecture was unsolved for 46 years.

\medskip
For Riemann integrable functions $f$
by a result of
H. Weyl \cite{[We]} from 1916,
there is a positive answer to the above question.

Raikoff \cite{[Raik]} (1936)
showed that if $f$ is
an arbitrary function
 with period one and integrable on $[0,1)$
then for an arbitrary
integer $q>1$
\begin{equation}\label{*rai}
\frac{1}{k}\sum_{n=1}^{k} f( q^{n}x)\to \int_{0}^{1}f(t)dt,
\end{equation}
holds for almost every
 $x$.
Raikoff also has shown that for any $\{ \lll_{k} \}$ increasing sequence of natural numbers for the above functions
$\frac{1}{k}\sum_{n=1}^{k} f( \lll_{n}x)$ converges to $ \int_{0}^{1}f(t)dt$ in the mean,
that is 
$$\int_{0}^{1}\Big | \Big (\frac{1}{k}\sum_{n=1}^{k} f( \lll_{n}x)\Big )- \int_{0}^{1}f(t)dt\Big |dx \to 0 \text{ as }k\to\oo.$$

Later  Riesz in  \cite{[Rieszerg]} observed that \eqref{*rai} is a consequence 
of the Birkhoff Ergodic Theorem.

In several papers Kac, Salem,  Zygmund, and Koksma  (\cite{[KSZ]}, \cite{[Koka]} and \cite{[Kokb]}) provided    positive answers to the Khinchin conjecture with additional assumptions on the
Fourier coefficients of $f$.

On the other hand, Erd\H os in \cite{[Erdstr]} showed in 1949
that there exist a sequence
$\lll_{n}\to\oo$ and a $1$-periodic function $f:\R\to \R$, 
$f\in L^{1}[0,1]$
for which
\begin{equation}\label{*erd}
\frac{1}{ k}\sum_{n=1}^{k} f( \lll_{n}x)\to
\int_{[0,1]}f
\end{equation}
 does not hold for
almost all $x$.
In fact, he gave a function $f$ which satisfies in addition to the above assumptions,  $\int_{[0,1]}f=0$, $\int_{[0,1]}f^{2}=1$
and for this function for a suitable sequence of integers $\lll_{n}\to\oo$ one has
$$\limsup_{k\to\oo } \frac{1}{ k}\sum_{n=1}^{k} f( \lll_{n}x)=+\oo$$
for $\mu$ a.e. $x.$

Finally, in 1969
Marstrand \cite{[M]}
showed that the Khinchin conjecture is not true.
In fact, in his example the set $E$, with $f=\chi_{E}$ can be selected to be an open set.
 He also gave further examples of sequences $\{ \lll_{n} \}$ for which
 \begin{equation}\label{*khii}
\frac{1}{k}\sum_{n=1}^{k} f(\lll_{n}x)\to  {\mu}(E)
\end{equation}
holds  for a.e. $x$, and sequences for which 
for certain sets \eqref{*khii} does not hold for a.e. $x$. For example, he showed that in this respect the sequence of primes is like the sequence of the integers
in that \eqref{*khii}  does not hold for certain open sets.

Marstrand's paper was continued by Baker in
\cite{[Baker]}. In this paper he was interested in 
properties of those open sets $E$ which can serve as an example for a sequence $\{ \lll_{n} \}$, for which  \eqref{*khii}
does not hold  for a.e. $x$.

Results of Marstrand and Baker were further refined by Nair in  four papers 
\cite{[Nair1]}, \cite{[Nair2]}, \cite{[Nair3]} and \cite{[Nair4]}.
For example in \cite{[Nair1]} it is shown that if the
 the sequence $\{ \lll_{n} \}$ is the increasing enumeration of a finitely
generated multiplicative subsemigroup of the positive integers, then for all $f \in L^1[0,1)$ which are $1$-periodic
\eqref{*erd} holds   for almost all $x$.

A different counterexample for the Khinchin conjecture  was given by J. Bourgain \cite{[Bou]} by using his entropy method.
His method permits one to unify certain counterexamples. Among other things the following question of Erd\H os
concerning a weaker version of the
Khintchine conjecture was also disproved: 
\begin{problem}\label{*bouerd}
{\it Is it true that given a
measurable set $E\sse \R$ with period $1$ for almost all $x$ the set $\{ n \in\N  : nx\in E \}$ has a
logarithmic density, that is letting $f(x)=\chi_{E}(x)$
$$\frac{1}{\log k}\sum_{n=1}^{k}\frac{f(nx)}n \to \mmm(E)?$$}
\end{problem}

Another counterexample to the Khinchin conjecture was given by Quas and Wierdl \cite{[QW]}.
In that paper several questions were asked. Here we quote just one of them (Question 7.4 of  \cite{[QW]}):
\begin{question}\label{*qqw}{\it
 Does there exist $f \in L ^1 ([0, 1))$, such that 
 $$\frac{1}{k \log k} \sum_{n=1}^{k}
 f (\{nx\})$$ diverges
almost everywhere?}
\end{question}
Quas and Wierdl also answered  a question of Nair 
\cite{[Nair1]}
 concerning the
Khinchine averages taken along a multiplicative subsemigroup of the natural numbers
rather than all natural numbers. 
They showed that the averages $$\frac{1}{\#(G \cap
[1, k ]))} \sum
_{\{n\in G : n\leq k \}} f (\{nx\})\to \int_{[0,1]}f$$ converge for a.e. $x$ for all $f \in L ^1[0,1)$ if and only if the
semigroup $G$ is a finitely generated subsemigroup. The "if" part of this result was proved by Nair
and the only if part is due to Quas and Wierdl.  

Functions 
\begin{equation}\label{*gap}
\text{$f:\R\to\R$, $1$-periodic  with $\int_{[0,1)}f=0$, $\int_{[0,1)}f^{2}<+\oo$}
\end{equation}
were considered by
Gapo\v{s}kin
\cite{[Gaplac]}, \cite{[Gapsys]} and  \cite{[Gapcert]}. 
He was interested in the asymptotic properties of the sequence $f(\lll_n x)$ for rapidly increasing sequences
$ \lll_{n}$ of integers. For example he proved that
for any $f$ satisfying \eqref{*gap} there exists a sequence $ \lll_{n} $
such that the sequence $f(\lll_{n } x)$ imitates the properties of independent random
variables. He  also gave a necessary and sufficient condition for $f(\lll_n x )$ to obey the central limit theorem. 
This line of research was followed by Berkes in \cite{[Bee]}.

Berkes in \cite{[Berk12]} and Berkes and Weber in 
 \cite{[BWa]} 
 study conditions under which the series
\begin{equation}\label{*BWa}
\sum_{n=1}^{\oo} c_{n}f(\lll_{n}x)
\end{equation}
defines an element of $L^2(\T)$ and converges to it in $L^2(\T)$, or pointwise for almost all $x$ in $\T$, where $\{\lll_n\}$ is an increasing sequence in $\N$, and $\{c_n\}\in \ell^2$. 
These conditions might be assumptions on $c_{n}$, or assumptions about $f$ like
the
H\"older exponent, properties of Fourier coefficients etc.
The nice survey \cite{[BWa]} is continued by many more recent articles, like
\cite{[Ab11]}, \cite{[W11]}, \cite{[Ai12]}, \cite{[BW14]}, \cite{[BWb]} \cite{[W15]}, \cite{[ABW]}, \cite{[ACBS]} and \cite{[W16]}.

Beck in \cite{[Beck]}   shows how to “save” Khinchin’s conjecture in the continuous case. It means that he uses instead of 
a sequence (arithmetic progression) a continuous
torus line. See also \cite{[BeckDon]}.

\section{$\sum f(nx)$, the non-periodic case}\label{*secfnx}

The Khinchin conjecture dealt with periodic functions $f$.
On the other hand Question \ref{*prWei} was about measurable functions defined 
on the positive semiaxis.

The following problem was asked by
 K. L. Chung in 1957 in the
 American Mathematical Monthly:
\medskip

{\bf Problem 4670 of the American Mathematical Monthly  \cite{[CHF]}:}
{\it If $f\in C[0,\oo)$, $f\geq 0$, and
$\int_{0}^{\oo}f(x)dx=\oo$ then
show that there exists $x>0$, for which
$\sum_{n=1}^{\oo}f(nx)=\oo.$}
\medskip

Fine and Hyde provided two solutions in \cite{[CHF]}.
They also remarked  that the answer to this problem also
provides a solution to
\medskip

{\bf Problem 4605 of the American Mathematical Monthly \cite{[AMM4605]}:}
{\it Given an open, unbounded set of positive reals. Prove that there exists a real number
such that infinitely many integral multiples of it lie in the set.}
\medskip

Indeed, if $E$ satisfies the assumptions of this problem then one can select 
 a non-negative continuous function $f$ such that it vanishes outside $E$ and
$\int_0^{\oo}f(x)dx = \oo$. For this function, $\sum_{n=1}^{\oo}f(nx)=\oo$ implies that $nx\in E$ for infinitely many $n$.

The above problems were followed by this one:

{\bf Problem 4727 of the American Mathematical Monthly \cite{[AMM4727]}:}
{\it Find a function $f(x),$ upper semi-continuous and
non-negative on $[0, \oo),$ bounded on each finite interval $(0, T),$ such
that $\int_{0}^{\oo}f(x)dx=\oo$ and $\sum_{n=1}^{\oo}f(nx)<\oo$ for every $x>0.$}

\medskip

These problems motivated 
Lekkerkerker in \cite{[L]} to work on these type of questions.
Supposing that $f:[0,\oo)\to[0,\oo)$ is a measurable function the following theorems hold:

\begin{theorem}[\cite{[L]}]\label{*thlek1}
If $\int_{0}^{\oo}f(x)dx<\oo$ then $\sum_{n=1}^{\oo}f(nx)<\oo$ for almost all
$x$.
\end{theorem}

\begin{theorem}[\cite{[L]}]\label{*thlek2}
If $f(x)$ is Riemann integrable over every finite interval and $\int_{0}^{\oo}f(x)dx=\oo$
then there exist $x>0$ such that $\sum_{n=1}^{\oo}f(nx)=\oo$.
\end{theorem}

In a footnote of \cite{[L]} it is mentioned that Theorem \ref{*thlek2} is due to  N. G. De Brujin.

\begin{theorem}[\cite{[L]}]\label{*thlek4}
 If $\mmm(E)=\oo$ and $E$ is Jordan measurable, then the set of
numbers $x > 0$ for which $nx\in E$ for infinitely many positive integers $n$
  is everywhere dense in the interval $[0,\oo)$
and has the power of the continuum.
\end{theorem}

\begin{theorem}[\cite{[L]}]\label{*thlek3}
 There exists a Jordan measurable set $E$ (consisting of
an infinite sequence of disjoint intervals) with $\mmm(E)=\oo$ such that 
$nx\in E$ for at most finitely many positive integers $n$
for almost all $x$ and, moreover, 
$nx\in E$
 for no positive integer $n$
for a set of $x$s of infinite measure.
\end{theorem}

In  \cite{[L]} it is also mentioned that the construction behind the example of this theorem is due to H. Kesten. It is also clear that by removing a set of zero Lebesgue measure 
from the set in Theorem \ref{*thlek3} one can obtain a Lebesgue measurable set 
for which the rest of the assumptions of Theorem \ref{*thlek4} hold, but the conclusion of this theorem does not.

Lekkerkerker in \cite{[LII]} considered higher dimensional generalizations.
In $\R^{d}$, $d\geq 2$. Instead of $nx$, $n\in \N$ he took
points of the form $n_{1}x^{1}+...+n_{d}x^{d}$, with $(n_{1},...,n_{d})\in \Z^{d}$ and
points $x^{j}\R^{d}$, $j=1,...,d$ which are independent in the sense, that they do not lie on
a $(d-1)$-dimensional hyperplane through the origin.

Sets considered  in Theorem \ref{*thlek3} motivated 
J. A. Haight \cite{[H1]}
to prove that there exists
a set
$E {\subset} (0, {\infty})$ of infinite measure, for which for all
 $x,y\in E,\  x\not=y$
the ratio $x/y$ is not an integer, and furthermore
{
 for all $x>0$,
$nx\not\in E$ if $n$ is sufficiently large.}

This implies that if
$f(x)=\chi_{E}(x)$,
the characteristic function of
 $E$ then
$\int_{0}^{{\infty}}f(x)dx= {\infty}$
and $\sum_{n=1}^{{\infty}}f(nx)< {\infty}$
everywhere.

In the introduction of \cite{[H2]} Haight remarked that sets of Lekkerker from Theorem \ref{*thlek3} were also constructed independently by W. M. Schmidt in  \cite{[Sch]}.
Schmidt in that paper asked   whether it was possible to have a set of infinite measure
with no two points having integral ratio. This question was answered independently by Haight in \cite{[H1]} and by Szemer\'edi \cite{[Szemer]}. 

In \cite{[H2]} the main result of \cite{[H1]} was generalized the following way:

{\it Consider any discrete enumerable set $G\sse (0,+\oo)$. A measurable set $E\sse (0,\oo)$ is constructed with the following properties:\\
 $\bullet$ If $x\in E,$ then $gx\not\in E$ for any $g\in G$;\\
 $\bullet$ if $y\in(0,\oo)$ then $yg\in E,$ $g\in G,$ has only a finite number of solutions $g$;\\
$\bullet$ the Lebesgue measure of $E$ is infinite. }
 
 The result of \cite{[H1]} was the special case $G=\N$.

Haight in  \cite{[H1]} and also in \cite{[H2]} stated a version of Problem \ref{*prWei}
in which $f$ is a characteristic function of a measurable set $E$.
Haight in \cite{[H2]} also generalized Problem \ref{*prWei} a bit in
the following way:
\begin{problem}\label{*prWeiH}
{\it 
Suppose that a discrete enumerable set $G\sse (0,+\oo)$ is given, $E$ is any Lebesgue measurable set of infinite measure and
$$N(x, E, G) = \sum_{\lll \in G}\chi_{E}(\lll x)= \# \{\lll : \lll x \in E, \lll \in G\}.$$
 Is it true that either $N(x, E,G) = \oo$ for almost all $x$
or that $N(x, E,G) < \oo$ for almost all $x$?}
\end{problem}

In  \cite{[BM1]} with D. Mauldin
we answered Question \ref{*prWei}, the Haight--Weizs\"aker problem.

\begin{theorem}\label{*HWth} There exists a
measurable function $f:(0,+ {\infty})\to \{0,1\}$ and
two nonempty intervals $I_{F}, \ I_{{\infty}} {\subset} [{1\over 2},1)$
such that for every $x\in I_{{\infty}}$ we have $\sum_{n=1}^{{\infty}}f(nx)=+ {\infty}$
and for almost every $x\in I_{F}$ we have $\sum_{n=1}^{{\infty}}f(nx)<+ {\infty}.$
The function $f$ is the characteristic function of an open set
$E$.\end{theorem}

The key tool in the proof of Theorem
\ref{*HWth}
was Kronecker's Theorem on simultaneous inhomogenous approximation
\cite{[C]}, p. 53. Here we state a special case of it which was used in
\cite{[BM1]}. ($||x||$ denotes the distance to the closest integer.)\medskip

\begin{theorem}\label{kronecker}
Assume $\theta_{1},...,\theta_{L}\in\R$ and
$(\alpha_{1},...,\alpha_{L})$ is a real vector. The following two statements
are equivalent:
\begin{itemize}
\item{A)} For every $\epsilon>0$, there exists $p\in\Z$ such that
$$||\theta_{j}p-\alpha_{j}||<\epsilon , \hbox{ for } 1\leq j\leq L.$$
\item{B)} If $(u_{1},...,u_{L})$ is a vector consisting of integers and
$$u_{1}\theta_{1}+...+u_{L}\theta_{L}\in\Z,$$
then $$u_{1}\alpha_{1}+...+u_{L}\alpha_{L}\in\Z.$$
\end{itemize}
\end{theorem}

\section{Translates of non-negative functions, type 1 and type 2 sets}\label{*sect12}

Since $\sum_{n=1}^{{\infty}}h(nx)=\sum_{n=1}^{{\infty}}h(e^{\log x+\log n})$,
then using the function $f=h\circ    \exp$ defined on $ {\ensuremath {\mathbb R}}$ and $\Lambda=\{\log n:
n=1,2,\ldots \}$  
we can rephrase Question \ref{*prWei} in the following way:
\begin{problem}\label{*prWeib}
{\it Suppose $f:\R \to [0,+\infty)$ is a measurable function.
Is it true that
$\sum_{{\lambda}\in {\Lambda}}f(x+ {\lambda})$
either converges (Lebesgue) almost everywhere or diverges almost everywhere, i.e.
is there  a zero-one law for $\sum f(x+\log n)$?}
\end{problem}
In the spirit of Question \ref{*prWeiH} one can ask the same question for other sets
$\LLL$ where 
\vskip-0.6cm
\begin{align}
&\text{\it $\LLL$ is a countably infinite  subset of $\R$  which has the unique}\label{*proplll}\\
&\text{\it accumulation point $+\oo$,} \nonumber
\end{align}
\vskip-0.2cm
\noindent of course elements of such a set $\LLL$ can be 
arranged as  $\{\lambda_1,\lambda_2,\ldots\}$,
a strictly increasing sequence tending to infinity.

 {\it In the rest of this paper if we write 
$\LLL$, $\LLL'$ or $\widetilde{\Lambda}$ property \eqref{*proplll} is always assumed.}

\medskip

In \cite{[BKM1]} Kahane suggested the definition of type 1 and type 2 sets.

\begin{definition}\label{def1} 
Given $\Lambda$ and
a measurable $f: {\ensuremath {\mathbb R}}\to [0,+ {\infty})$  let
$$ s_f(x)  = s(x)=\sum_{\lambda\in\Lambda}f(x+\lambda).$$
Consider the complementary subsets of $ {\ensuremath {\mathbb R}}$:
\begin{equation}\label{*CD}
C=C(f,\Lambda)=\{x: s(x)< {\infty}\},\text{ and }
D=D(f,\Lambda)=\{x:s(x)= {\infty}\}.
\end{equation}
The set $\Lambda$ is type $1$ if, for every $f$, either
$C(f,\Lambda)= {\ensuremath {\mathbb R}}$ a.e. or $C(f,\Lambda)= {\emptyset}$ a.e.~(or equivalently
$D(f,\Lambda)= {\emptyset}$ a.e. or $D(f,\Lambda)= {\ensuremath {\mathbb R}}$ a.e.). Otherwise, $\Lambda$
is type $2$.
For type $2$ sets there are non-negative measurable {\it witness functions}
$f$ such that both $C(f,\Lambda)$ and $D(f,\Lambda)$ are of positive measure.
\end{definition}

That is, for type 1 sets there is a ``zero-one'' law for the almost everywhere
convergence properties of the series $\sum_{\lambda\in\Lambda}f(x+\lambda)$,
while for type 2 sets there is no such law.

Although we will see several necessary and several sufficient conditions
in this survey article, the main unsolved problem in this area is:
\begin{question}\label{*qtype12}
{\it
Give a (reasonable) necessary and sufficient condition for a set $\LLL$ being type $2$, that is characterize type $2$ sets.}
\end{question}

We denote by $C^{+}_{0}( {\ensuremath {\mathbb R}})$ the non-negative continuous functions on $ {\ensuremath {\mathbb R}}$  tending to zero in $+ {\infty}$. 

 We recall here 
from \cite{[BKM1]}
the theorem concerning the Haight--Weizs\"aker problem. This contains the additive version of the result of Theorem \ref{*HWth}
with some extra information.

Using this terminology in \cite{[BKM1]} we stated the following version of
Theorem \ref{*HWth}:
\begin{theorem}\label{*BKMHWth} {The set $\Lambda=\{\log n:
n=1,2,... \}$ has type $2$. Moreover, for some $f\in
C_{0}^{+}( {\ensuremath {\mathbb R}}),$ $C(f, {\Lambda})$ has full measure on the half-line
$(0,\infty)$ and $D(f, {\Lambda})$ contains the half-line $(-\infty,0)$.
If for each $c, \int_c^{+\infty}e^yg(y)dy < +\infty$, then $C(g,\Lambda) =  {\ensuremath {\mathbb R}}$ a.e.
If $g\in C_{0}^{+}( {\ensuremath {\mathbb R}})$ and $C(g,\Lambda)$ is not of the first (Baire) category, then
$C(g, {\Lambda})= {\ensuremath {\mathbb R}}$ a.e. Finally, there is some $g\in C_{0}^+( {\ensuremath {\mathbb R}})$ such that $C(g,\Lambda) =  {\ensuremath {\mathbb R}}$ a.e. and
$\int_0^{+\infty}e^yg(y)dy = +\infty$.
}
\end{theorem}
\begin{center}
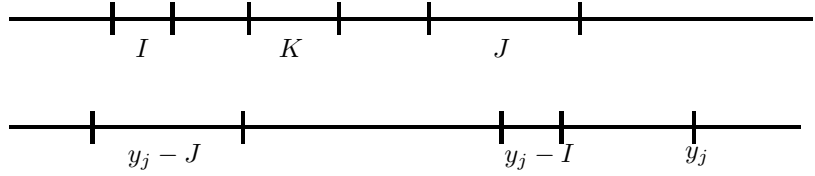
\begin{figure}[!t]
\hskip-1.6cm
\setlength{\unitlength}{0.8mm}
\linethickness{1.25pt}
\begin{picture}(153.67,25)
\newsavebox{\foldera}
\savebox{\foldera}
  (0,0)[bl]{
\put(20.33,18.67){\line(1,0){133.33}}
\put(37.33,21.33){\line(0,-1){5.33}}
\put(47.33,21.33){\line(0,-1){5.33}}
\put(60.00,21.33){\line(0,-1){5.33}}
\put(75.00,21.33){\line(0,-1){5.33}}
\put(90.00,21.33){\line(0,-1){5.33}}
\put(115.00,21.33){\line(0,-1){5.33}}
\put(102.00,13.78){\makebox(0,0)[cc]{$J$}}
\put(67.00,13.78){\makebox(0,0)[cc]{$K$}}
\put(42.33,13.78){\makebox(0,0)[cc]{$I$}}
}

\newsavebox{\folderb}
\savebox{\folderb}
  (0,0)[l]{
 \put(18.33,20.00){\line(1,0){131.34}}
\put(132.00,22.67){\line(0,-1){5.34}}
\put(132.50,15.11){\makebox(0,0)[cc]{$y_j$}}
\put(110.00,22.67){\line(0,-1){5.34}}
\put(100.00,22.67){\line(0,-1){5.34}}
\put(57.00,22.67){\line(0,-1){5.34}}
\put(32.00,22.67){\line(0,-1){5.34}}
\put(44.00,15.11){\makebox(0,0)[cc]{$y_j-J$}}
\put(106.33,15.11){\makebox(0,0)[cc]{$y_j-I$}}  
}

\put(0,8){\usebox{\foldera}}
\put(2,0){\usebox{\folderb}}
\end{picture}
\caption{The intervals in Theorem \ref{*BKMth5}.}
\label{*ffig1}
\end{figure}
\end{center}
For sets $A,B \subset  {\mathbb R}$ we define the Minkowski sum
$A+B=\{a+b: \  a\in A,\  \, b \in B\}$ and $A-B=\{a-b: \  a\in A, \  \, b \in B\}$.

In \cite{[BKM1]} 
and in several other results in later papers 
the key ingredient in the proofs was the following:
\begin{theorem}[Theorem 5 of \cite{[BKM1]}] \label{*BKMth5}
Suppose that there exist three intervals $I$, $J$, $K$ such that $J=K+I-I$,
$I$ is to the left of $J$ (see 
Figure \ref{*ffig1}), and $\hbox{dist}(I,J)\geq \mu(I),$ and
two sequences $(y_{j})$
and $(N_{j})$ tending to infinity ($y_{j}\in \R^{+}$,
$N_{j}\in\N$) such that, for each $j$, $y_{j}-I$
contains a set of $N_{j}$ points of $\Lambda$ 
independent from $\Lambda \cap (y_j-J)$ in the sense that
the additive groups generated by these sets
have only 0 in common. Then $\Lambda$ has type
2. 
Moreover, for some $f\in C_{0}^{+}(\R)$, $D(f,\Lambda)$ contains $I$ and
$C(f,\Lambda)$ has full measure on $K$.
\end{theorem}

For the proof of this theorem instead of Theorem \ref{kronecker} it is better to use Kronecker's theorem 
in this "homogeneous harmonic analysis" version from Katznelson's book (see Theorem 9.1 in Chapter VI of \cite{Katznelson}):
\begin{theorem}\label{KronKatz}
Let $\lll_1,\ldots,\lll_n$ be real numbers, independent
over the rationals. 
Let $\alpha_1,\ldots,\alpha_n$ be real numbers and $\eee>0$.
Then there exists a real number $x$ such that 
$$
\left|e^{i\lll_j x} - e^{i\aaa_j}\right|<\eee, \qquad  j=1,\ldots,n.
$$
\end{theorem}

In \cite{[BKM1]} we showed that type $2$ is typical in the Baire category sense.
This matches the general opinion that typical objects are "bad", like typical continuous functions are nowhere differentiable.

To speak about Baire category we need to define 
 a topology $ {{ \cal T}}$
on the space  $ {{ \cal D}}$ of infinite discrete subsets of $ {{\R}}^+$.
Let $\Lambda \in  {{ \cal D}}$ and let $r_{n}\in  {{\R}}^+$ for $n=1,2,...$.
The set $\Lambda'$ belongs to the
$ {{ \cal N}}((  \ell
_{1},r_{1}),(  \ell _{2},r_{2}),...)$
neighborhood of
$\Lambda$ if we can order the elements of $\Lambda$ into a sequence
$\{  \ell _{1},  \ell _{2},...\}$ and we can order the elements of $\Lambda'$ into a sequence
$(  \ell '_{n})$ such that $  \ell' _{n}\in (  \ell _{n}-r_{n},  \ell
_{n}+r_{n}).$ The topology
$ {{ \cal T}}$
will be generated by these neighborhoods. 
The space $( {{ \cal D}}, {{ \cal T}})$ is a Baire
space (see for example \cite{[BKM1]}).

\begin{theorem}[\cite{[BKM1]}]\label{*BKMth6}
 The sets of type $2$ form a dense
open subset in $ {{ \cal T}}$, while sets of type $1$ form a closed
nowhere dense set. Therefore type $2$ is typical in the Baire
category sense.
\end{theorem}

\section{Dyadic examples}\label{*sectdyad}

There is a class of sets $\LLL$ for which we have a complete characterization of type $1$ and $2$ sets.

We start with two examples from \cite{[BKM1]}:

\begin{example}\label{*exdyada}
Set $\Lambda=\cup_{k \in {{{\ensuremath {\mathbb N}}}}}\Lambda_{k},$ where $\Lambda_{k}=2^{-
k} {{{\ensuremath {\mathbb N}}}}\cap [k,{k+1}).$
In Theorem 1 of \cite{[BKM1]} it is proved that $\Lambda$ is type $1$.
In fact, in a slightly more general version it is shown that
if $(n_{k})$ is a strictly increasing sequence of
positive integers  and $\Lambda=\cup_{k\in  {{{\ensuremath {\mathbb N}}}}}\Lambda_{k}$ where $\Lambda_{k}=2^{-
k} {{{\ensuremath {\mathbb N}}}}\cap [n_{k},n_{k+1})$
then $\Lambda$ is type $1$.
\end{example}

\begin{example}\label{*exdyadb}
Let
$(n_{k})$ be a given increasing sequence of positive integers.
 By Theorem 3 of \cite{[BKM1]}
there is
an increasing sequence of integers $(m(k))$ such that the set
$\Lambda=\cup_{k\in  {{{\ensuremath {\mathbb N}}}}}\Lambda_{k}$ with $\Lambda_{k}=2^{-m(
k)} {{{\ensuremath {\mathbb N}}}}\cap [n_{k},n_{k+1})$ is type $2$.
\end{example}

Theorem \ref{mkthm} from \cite{[BHMVa]} is a much sharper version of Example \ref{*exdyadb} since it gives a
necessary and sufficient condition for a set obtained by the ``dyadic" construction
being type $1$ (or type $2$). It also shows that for type $2$ a faster rate of increase of elements in the "next dyadic block" is necessary.

\begin{theorem}\label{mkthm}
Suppose that $(m_k)$ and $(n_k)$ are strictly increasing sequences of positive integers.  For each $k \in  {\mathbb {N}}$, define $\Lambda_k=2^{-m_k} {\mathbb {N}} \cap [n_k,n_{k+1})$ and let $\Lambda=\cup_{k=1}^\infty \Lambda_k$.  Define $  \displaystyle  M=\sup_k\{m_{k+1}-m_k\}$.  Then $\Lambda$ is type 1 if and only if $M < \infty$.
\end{theorem}

In \cite{[BHMVr]} we studied the effect of randomly deleting elements
of $ {\Lambda}$.
Let $0<p<1$.  Then we say that $\Lambda\subset \widetilde{\Lambda}$ is chosen with probability $p$ from $\widetilde{\Lambda}$ if for each $\lambda \in \widetilde{\Lambda}$ 
independently the probability that $\lambda \in \Lambda$ is $p$.

It is clear that almost surely if $\Lambda\subset \widetilde{\Lambda}$ is chosen with probability $p$ from $\widetilde{\Lambda}$ then $ {\Lambda}$ 
satisfies \eqref{*proplll}.

Let
$
\widetilde{\Lambda}=\bigcup_{k=1}^\infty (2^{-k} {\mathbb {N}} \cap[k,k+1)).
$
 We know from
Example \ref{*exdyada} that $ {{\widetilde {\Lambda}}}$ is type $1$.
By Theorem \ref{prob p theorem}
if $\Lambda$ is chosen with probability $p$ from $\widetilde{\Lambda}$
then almost surely $\Lambda$ is type $1$. This again matches the heuristic idea that
more elements can cause trouble, deleting random elements from a good set keeps it good.

\begin{theorem}[Theorem 4.3 of \cite{[BHMVr]}]\label{prob p theorem}
Suppose that $0<p<1$ and $\Lambda$ is chosen with probability $p$ from $\widetilde{\Lambda}=\cup_{k=1}^\infty (2^{-k} {\mathbb {N}} \cap[k,k+1))$.  Then $\Lambda$ is type $1$
almost surely.
\end{theorem}

The next theorem might be surprising, since it contradicts the previous heuristic idea.  It may happen that
type $1$ sets are converted into type $2$ sets by random deletion:
\begin{theorem}[Theorem 4.5 of \cite{[BHMVr]}]\label{*randt2} Suppose that $(m_k)$ and $(n_k)$ are strictly increasing sequences of
positive integers. For each $k\in\mathbb{N}$, define $\Lambda_k=2^{-m_k} \mathbb{N} \cap [n_k, n_{k+1})$ and let $ {{\widetilde {\Lambda}}}=\bigcup_{k=1}^{\infty}\Lambda_k$. Moreover, fix $0<p<1$ and suppose that $\Lambda$ is chosen  with probability $p$ from $ {{\widetilde {\Lambda}}}$.
Set $q=1-p$.
 For fixed $(m_k)$, if $(n_k)$ tends to infinity sufficiently fast
 then  almost surely $ {\Lambda}$ is type $2$. Notably, if the series 
 \begin{equation}\label{*nkgr}
 \sum_{k=1}^{\infty}1-\left(1-{q}^{2^{m_k}}\right)^{n_{k+1}-n_k} \text { diverges}
\end{equation}
  then  almost surely $ {\Lambda}$ is type $2$. \end{theorem}
 
By Theorem \ref{mkthm} if  $M=\sup_k\{m_{k+1}-m_k\}<+\oo$ we obtain a type $1$ set.
Next we can choose $n_{k}$ such that \eqref{*nkgr} is satisfied. This way
using Theorem \ref{*randt2} by random deletion we obtain a type $2$ set from a type $1$ set.

 \medskip
 
In \cite{[BHMVa]} another variant of modifying dyadic constructions was considered.
This time the result is in line with the heuristic idea that adding elements, especially algebraically independent elements, results in bad, type $2$ sets.

\begin{theorem}\label{type1plusalg}
Let $\{m_k\}_{k=1}^N$ and $\{n_k\}_{k=1}^N$ be strictly increasing sequences of positive integers, where either $N \in  {\mathbb {N}}$,
 or $N=\infty$.  If $N \in  {\mathbb {N}}$ we define $n_{N+1}=\infty$. Define
$$\Lambda_1=\cup_{k=1}^N (2^{-m_k} {\mathbb {N}} \cap [n_k,n_{k+1})).$$
Let $\Lambda_2=\{\alpha_k:\ k\in {\mathbb {N}}\}$ be a set of irrational numbers
 independent over the rationals, where $\alpha_k \nearrow \infty$.
Then $\Lambda_*=\Lambda_1 \cup \Lambda_2$ is type 2.
\end{theorem}

Observe that $ {\Lambda}_{1}$ can be any of the sets from Examples \ref{*exdyada}
or \ref{*exdyadb}, hence $ {\Lambda}_{1}$ can be a type $1$ set which is converted in this case
into a type $2$ set after we add the independent numbers.
Theorem \ref{*BKMth5} also suggests that lots of independence 
implies type $2$. However this is false, as the next theorem illustrates:

\begin{theorem}[Theorem 3.4 of \cite{[BHMVa]}]\label{alg indep type 1}
There exists a set $\Lambda$ which is type 1 and which includes infinitely many numbers independent over the rationals.
\end{theorem}

One might believe that for type 2 sets $\Lambda$ the sets $C=C(f, {\Lambda})$, or $D=D(f, {\Lambda})$ are always half-lines
if they differ from $ {\ensuremath {\mathbb R}}$.
Indeed in \cite{[BKM1]} we obtained results in this direction.
A number $t>0$ is called a translator of $ {\Lambda}$ if $( {\Lambda}+t) {\setminus}
 {\Lambda}$ is finite.

 \medskip
 
 {\it Condition $(*)$ is said to be satisfied\\
 if $T( {\Lambda})$, the
countable additive
semigroup of translators of $ {\Lambda}$, is dense in $ {\ensuremath {\mathbb R}}^{+}$.}

\medskip

We recall:

\begin{proposition}[Proposition 3 of \cite{[BKM1]}]\label{*proptrans}
Suppose that condition $(*)$ is
satisfied ($\Lambda$ has arbitrarily small translators).
Then the topological closure of $C$ (resp. $D$) is either $ {\emptyset}$, or
$ {\ensuremath {\mathbb R}}$, or else a closed right half-line (resp. left half-line). The same
holds for the support of $\pmb{1}_{C}$ (resp. $\pmb{1}_{D}$) meaning the
smallest closed set $\pmb{C}$ carrying $C$ (resp. $\pmb{D}$ carrying $D$) except for a null set.
The interior of $\pmb{C}$ (resp. $\pmb{D}$) is either $ {\emptyset}$, or $ {\ensuremath {\mathbb R}}$, or else an
open right (resp. left) half-line.
\end{proposition}

 From the sufficient condition of Theorem \ref{*BKMth5} one might think that for
sets satisfying condition($*$) some sort of independence
determines whether the set is type $2$. However, Example \ref{*exdyadb} shows that this
is not the case. There are type $2$ sets $\Lambda$ which satisfy condition($*$) and for which any two elements of $\Lambda$ 
are dependent over the rationals.

 \medskip

In Theorem \ref{th1} we can find the first example where $C(f, {\Lambda})$ does not equal
$ {\emptyset}$, $ {\ensuremath {\mathbb R}}$, or a left half-line modulo sets of measure zero.

In \cite{[BM2]} another variant of the dyadic construction was considered.
Suppose that $ {\alpha}\in(0,1)$, and let
$ {\Lambda}^{{\alpha}^{k}}:=\cup_{k=1}^{{\infty}} {\Lambda}_{k}^{{\alpha}^{k}},$
$ {\Lambda}_{k}^{{\alpha}^{k}}
= {\alpha}^{k} {\ensuremath {\mathbb Z}}\cap [n_k,n_{k+1})$.

If $\alpha\not\in \Q$ then using
Theorem \ref{*BKMth5} one can show that
${\Lambda}^{{\alpha}^{k}}$
is of type $2$.

If $ {\alpha}=\frac{1}q$ for some $q\in \{2,3,...  \}$, 
then a slight modification of the proof of
Theorem 1 of \cite{[BKM1]} (see Example \ref{*exdyada} in this paper) shows that  $ {\Lambda}^{(\frac{1}q)^{k}}$ is of type $1$ and condition $(*)$ is satisfied.

Interestingly the case when $ {\alpha}=\frac{p}q$ with
$(p,q)=1,$ $p,q>1$, $p<q$ is different, these rational numbers are behaving like the irrationals the set $\LLL= {\Lambda}^{{(\frac {p} {q} )}^{k}}$ is of type $2$.
The precise result is the following:
 \begin{theorem}[\cite{[BM2]}]\label{th1}
Assume the integer $r>8$, $(r,p)=1,$ and $(r,q)=1$.
 The set $ {\Lambda}$ defined above is of type 2. Moreover, there
exists
a characteristic function $f: {\ensuremath {\mathbb R}}\to [0,\oo)$ such that almost every
point in $[ {\frac {1} {r}}, {\frac {2} {r}}]$ belongs
to $C(f, {\Lambda})$,$$\text{ $\mu\Big (D(f, {\Lambda})\cap [1+ {\frac {1} {r}},1+ {\frac {2} {r}}]\Big )>\frac{1}{4r}$
and $\mu \Big (D(f, {\Lambda})\cap [-1+ {\frac {1} {r}},-1+ {\frac {2} {r}}]\Big )>\frac{1}{4r}.$}$$
\end{theorem}

So the rationals $ {\alpha}=\frac{p}q$, with $p>1$  behave like irrationals.

 Even in the special case when $p=2$,
$q=3$ the proof of Theorem \ref{th1} is not obvious.
For a while it seemed that we
needed some information on the distribution of $ \{
(\frac{3}2)^{k}  \}$,
(where $ \{\cdot   \}$ denotes the fractional part).
At that time we did not have information about the background of this question.
Finally on Mathscinet we found references to Choquet's papers
  \cite{[Choq80]}, \cite{[Choq80b]} that in 1980 it was not even known, whether $
\{(\frac{3}2)^{k}
  \}$ is uniformly distributed, or even dense in $[0,1]$. 
This question was extensively studied since then.
See, for example the monograph of Bugeaud \cite{[Bugeaud]}.
Here we recall some problems/conjectures from this book which seem to be still unsolved.

\begin{question}[Mahler \cite{[Mahler]}]\label{*qMahler}
{\it  There are no real numbers $\xi$ such that $0 \leq \{\xi(3/2)^n \} <
1/2$ for every positive integer $n.$}
\end{question}

\begin{question}\label{*qhkn}
{\it The sequence $(\{(3/2)^n\} )_{n\geq1}$ is dense modulo one.}
\end{question}

This has already been mentioned. Of course, the same problem with uniform density is even more difficult.
There is a weaker version of Problem \ref{*qhkn} which still seems to be unsolved:

\begin{question}[Mend\`es France \cite{[MendFr]}]\label{*qmendfr}
{\it The sequence $(\{(3/2)\}^n )_{n\geq1}$ has an irrational limit
point. }
\end{question}

For some recent progress related to the above questions we refer to 
\cite{[Dub19]} and
\cite{[Stru21]}.

\section{The structure of witness functions/sets for type $2$}\label{*sectwit}

The next question is about the structure of witness functions for type $2$ sets.
\begin{problem}[QUESTION 1 in \cite{[BKM1]}]\label{*q1bkm1}
{\it
Is it true that $\Lambda$ is type $2$ if and only if there
is a $\{0,1\}$ valued measurable function $f$ such that both
$C(f,\Lambda)$ and $D(f,\Lambda)$ have positive Lebesgue measure?}
\end{problem}

We gave the following answer to this question in \cite{[BHMVr]}:
\begin{theorem}\label{*thwit} Suppose that $\Lambda$ is type $2$, that is there exists a measurable witness function $f$ such that both $D(f,\Lambda)$ and $C(f,\Lambda)$ have positive measure. Then there exists a witness function $g$ which is the characteristic function of an open set and both $D(g,\Lambda)$ and $C(g,\Lambda)$ have positive measure.\end{theorem}

This definition is also from \cite{[BKM1]}:
\begin{definition}\label{asydens} The  set $\Lambda=\{{\lambda}_{1}, {\lambda}_{2},... \}$, $ {\lambda}_{1}< {\lambda}_{2}<...$ is asymptotically dense if $d_{n}= {\lambda}_{n}- {\lambda}_{n-1}\to 0$, or equivalently:
$$\forall a>0,\quad \lim_{x\to\infty}\#(\Lambda\cap [x,x+a])=\infty.$$
If $d_{n}$ tends to zero monotonically, we speak about
decreasing gap asymptotically dense sets.

If  $\Lambda$ is not asymptotically dense we say that it is asymptotically lacunary.
\end{definition}

\medskip

We note that the above notion of asymptotically lacunary is not the same as
the usual one of a ``lacunary sequence.''

Next we turn to questions related to the structure of convergence and divergence sets for
type $2$ $\LLL$s.
\begin{problem}[QUESTION 2 in \cite{[BKM1]}]\label{*q1bkm1} 
{\it Given open sets $G_1$ and $G_2$ when is it possible
to find $\Lambda$ and $f$ such that $C(f,\Lambda)$ contains $G_1$ and
$D(f,\Lambda)$ contains $G_2$?}
\end{problem}

In \cite{[BMV]} we gave an almost complete answer to this question.
\begin{theorem}[\cite{[BMV]}]\label{*thdflg}
There is a strictly monotone increasing unbounded sequence $(\lambda_0,\lambda_1,\ldots)=\Lambda$ in $ {\mathbb {R}}$ such that  $\lambda_{n}-\lambda_{n-1}$ tends to $0$ monotone decreasingly,
that is $ {\Lambda}$ is a decreasing gap asymptotically dense set,
such that for every open set $G\subset {\mathbb {R}}$ there is a function $f_G: {\ensuremath {\mathbb R}}\to [0,+ {\infty})$ for which
\begin{equation}\label{tetelbeli}
\mu\left(\left\{x\notin G :  \sum_{n=0}^\infty  f_G(x+\lambda_n)=\infty\right\}\right)=0, \text{   and   }
 \end{equation}
\begin{equation}\label{*conv}
 \text{   $\sum_{n=0   }^\infty f_G(x+\lambda_n) =\infty$ for every $x\in G$,}
\end{equation} 
moreover
$f_G=\chi_{U_G}$ for a closed set $U_G\subset {\mathbb {R}}$.
By \eqref{tetelbeli} and \eqref{*conv} we have $D(f_G, {\Lambda})\supset G$, and 
$C(f_G, {\Lambda})= {\ensuremath {\mathbb R}} {\setminus} G$ modulo sets
of measure zero. 
One can also select a $g_G\in C_{0}^{+}( {\ensuremath {\mathbb R}})$ satisfying \eqref{tetelbeli} and
\eqref{*conv} instead of $f_G$. 
\end{theorem}

Observe that in the above theorem we construct a universal $ {\Lambda}$ and for this
set, depending on our choice of $G$ we can select a suitable $f_{G}$ such that
$D(f_{G}, {\Lambda})=G$ modulo sets of measure zero.

This also leads to the following question:
\begin{question}\label{*typconvdiv}{\it
We have seen in Theorem \ref{*BKMth6} that the typical $\LLL$ is type $2$.
What can be said about the possible $C(f,\Lambda)$ and $D(f,\Lambda)$
sets for  typical $\LLL$s?  Are the typical $\LLL$s universal, like the ones
in Theorem \ref{*thdflg}?}
\end{question}

One can also state a problem about more universal $\LLL$s.
\begin{question}\label{*vert2}{\it
Characterize (or give sufficient conditions, and necessary conditions) of $\Lambda$s for which for every $H\subset\R$ there is an $f$ such that $C(f,\Lambda)=H$. (In this problem $H$ and $f$ are not necessarily measurable. Of course, one can state this problem with the assumption of measurability as well.)
}
\end{question}

Next we are interested in the possibility of saying something without an exceptional set 
of measure zero.
\begin{theorem}[\cite{[BMV]}]\label{*thdivG}
There exists an asymptotically dense infinite discrete set $ {\Lambda}$
 such that for any open set $G {\subset}  {\ensuremath {\mathbb R}}$ one can select an
 $f_{G}\in C^{+}_{0}( {\ensuremath {\mathbb R}})$  such that $D(f_G, {\Lambda})=G.$
\end{theorem}

In this theorem $D(f, {\Lambda})=G$, that is,  there is no exceptional set of measure
zero on which we do not know what happens. This also implies that if the interior
of $ {\ensuremath {\mathbb R}} {\setminus} G$ is non-empty then $C(f, {\Lambda})$ contains intervals.

It was remarked in \cite{[BKM1]}  that if the counting function
of $\Lambda,\  n(x)=\#\{{\Lambda}\cap [0,x]\}$ satisfies a condition of the type
$$
 \forall   \ell < 0,\   \forall a\in  {{{\ensuremath {\mathbb R}}}}, \ \ \limsup_{x \to
\infty}{n(x+  \ell+a) - n(x+a) \over n(x+  \ell) - n(x)} < +\infty  
$$
(as is the case for $\Lambda = \{\log n\}$) with witness functions $f\in C_{0}^{+}( {\ensuremath {\mathbb R}})$ then either $C(f,\Lambda)$ has full
measure on $ {{{\ensuremath {\mathbb R}}}}$ or $C(f,\Lambda)$ does not contain any interval.

It is possible to construct
$f \in C_0^+( {{{\ensuremath {\mathbb R}}}})$ such that both $C(f,\Lambda)$ and $D(f,\Lambda)$ have
interior points:
\begin{theorem}[\cite{[BKM1]}]\label{*BKMth4}
If $\Lambda$ is asymptotically lacunary, then $\Lambda$ is
type $2$. Moreover, for some $f\in C_{0}^{+}( {\ensuremath {\mathbb R}})$, there exist
intervals $I$ and $J$, $I$ to the left of $J$, such that $C(f,\Lambda)$
contains $I$ and $D(f,\Lambda)$ contains $J.$
\end{theorem}

The next example shows that asymptotical lacunarity of $\LLL$ is not needed for the previous theorem.  One can define decreasing gap asymptotically dense
$ {\Lambda}$s for which one can find nonnegative continuous $f$s such that both
$C(f, {\Lambda})$ and $D(f, {\Lambda})$  have interior points.

\begin{theorem}[\cite{[BMV]}]\label{*excint}
There exists a decreasing gap asymptotically dense $ {\Lambda}$ and an $f\in C^{+}_{0}( {\ensuremath {\mathbb R}})$ such that
$I_{1}=[1,2] {\subset} D(f, {\Lambda})$ and $I_{2}=[4,5] {\subset} C(f, {\Lambda})$.
\end{theorem}

Observe that in the above construction $I_{1} {\subset} D(f, {\Lambda})$ was to the left of $I_{2} {\subset} C(f, {\Lambda})$.
 The next theorem shows that for decreasing gap asymptotically dense $ {\Lambda}$s and continuous witness functions this situation cannot be improved. 
If $x$ is an interior point of $C(f, {\Lambda})$ then the half-line $[x, {\infty})$ intersects $D(f, {\Lambda})$ in a set of measure zero. 
As Theorem \ref{*thdivG} shows if we do not assume that $ {\Lambda}$ is of decreasing gap then
it is possible that $D(f, {\Lambda})$ has  a part of positive measure, even to the 
right of the  interior points of $C(f, {\Lambda})$. 

\begin{theorem}[\cite{[BMV]}]\label{*thcintl} Let $\Lambda$ be a decreasing gap and asymptotically dense set, and let $f: {\ensuremath {\mathbb R}}\to[0,+ {\infty})$ be continuous. Then if $x$ is an interior point of $C(f,\Lambda)$ then 
\begin{equation}\label{*zhalf}
 {\mu}\Big ([x,+ {\infty})\cap D(f, {\Lambda})\Big )=0.
\end{equation}
 \end{theorem}

This theorem
 implies that it is impossible to find such a universal $ {\Lambda}$, like in Theorem
\ref{*thdflg} with decreasing gaps.

In  \cite{[BKM1]} there were some simple observations concerning the dependence of
$C(f,\Lambda)$ and $D(f,\Lambda)$ on modifications of $f$ for a given $\Lambda$.

For example $C(f,\Lambda)=C(\min\{f,1\},\Lambda)$, so one can always assume that
$f$ is bounded.

One can additionally suppose that $f(x)$ tends to
$0$ at infinity
provided we are interested in $C(f,\Lambda)$ and $D(f,\Lambda)$ up to a nullset.

In \cite{[BKM1]} there was also a simple argument
based on the Borel--Cantelli lemma about modifications of witness functions 
with changes only modulo sets of zero measure in the convergence/divergence sets.
We state it here as the next lemma.
\begin{lemma}\label{*lemmodBKM}
Suppose that $ {\Lambda}$ is  type $2$ and $f$ is a bounded witness function for $ {\Lambda}$.
 If we modify $f$ on a set $E$ such that $\mu(E\cap(x,\infty))\leq{\epsilon(x)}$ where $\epsilon(x)$ is a positive decreasing function tending to $0$ at infinity, and satisfying
\begin{equation}\label{*modeq}
\sum_{l\in\mathbb{N}}\epsilon(l-K)\#(\Lambda\cap[l,l+1))<\infty,
\end{equation}
then the convergence and divergence sets in $[-K,K]$  for the modified function $\widetilde{f}$ do not change apart from a set of measure $0.$
\end{lemma}

Using the previous observations and this lemma one can obtain:
\begin{proposition}[Proposition 1 of  \cite{[BKM1]}]\label{*propco}
Given a Lebesgue measurable $f\geq 0$  
there exists $\widetilde{f}\geq 0$, continuous and tending to $0$ at infinity
($\widetilde{f}\in C_{0}^{+}( {{\R}})$) such that
$$C(f,\Lambda)=C(\widetilde{f},\Lambda)\hbox{ a.e., }D(f,\Lambda)=D(\widetilde{f},\Lambda)\hbox{ a.e.}$$
\end{proposition}

During the modification of witness functions quite often the Borel-Cantelli lemma is used and we end up with modified functions with the same convergence and divergence sets modulo sets of measure zero.
This motivates the following question about possible sharpening of Theorem \ref{*thdflg}.
\begin{question}\label{*qmodif}
{\it Suppose that $f$ is a measurable witness function with convergence set $C(f,\Lambda)$ for a type $2$ set $\LLL$. Is there always a characteristic function $g$ for which $C(g,\Lambda)=C(f,\Lambda)$?
}
\end{question}

A probably easier version is related to Proposition \ref{*propco}:
\begin{question}\label{*qmodif}
{\it Given $f\geq 0$ and Lebesgue measurable,
can we find $\widetilde{f}\geq 0$, tending to $0$ at infinity
such that
$$C(f,\Lambda)=C(\widetilde{f},\Lambda)\hbox{ and }D(f,\Lambda)=D(\widetilde{f},\Lambda)?$$
}
\end{question}
It is clear that in this question $\widetilde{f}$ cannot always be in $ C_{0}^{+}( {{\R}})$, since for arbitrary measurable functions $C(f,\Lambda)$ might be at a higher level of Borel hierarchy than the levels which can be reached by 
$C(f,\Lambda)$ for continuous functions. Of course, there is a further question:
\begin{question}\label{*hiera}
{\it If $C(f,\Lambda)$ is at a certain level of the Borel hierarchy for a measurable witness function $f$  what is the lowest level in the Baire classes from which we can always find
a witness function $\widetilde{f}$ such that  $C(f,\Lambda)=C(\widetilde{f},\Lambda)$?}
\end{question} 

\section{Type $1$ and $2$ sub and supersets}\label{*secsubsup}

The following theorem is a generalization of Theorem  \ref{*BKMth4}.
This again supports the heuristic idea that adding many points
leads to "bad" type $2$ sets. 

\begin{theorem}[\cite{[BHMVa]}]\label{ez is 2-es}
Let $ {\varepsilon}$ be a positive number. 
For every $n\in {\ensuremath {\mathbb Z}}$ we denote the cardinality of $\Lambda\cap [n {\varepsilon},(n+1) {\varepsilon})$ by $a_n$. If
\begin{equation}\label{gyors}
\limsup_{n\rightarrow\infty} \frac{a_n}{a_{n-1}} = \infty
\end{equation}
(where $\frac{0}{0}=0$ and $\frac{c}{0}=\infty$ if $c>0$), then $\Lambda$ is type 2.
\end{theorem}

In the next corollary we obtain a super "bad", super type $2$ $\LLL$. This means that
any superset of  $\LLL$ is also type $2$. The assumption in \eqref{fast}
is also in line with the heuristic idea that many points/fastly increasing number of points in $\LLL$ yield (super) type $2$ sets.
 
\begin{corollary}[\cite{[BHMVa]}]\label{suru 2-es}
For every $n\in {\ensuremath {\mathbb Z}}$ we denote the cardinality of $\Lambda\cap [n,n+1)$ by $a_n$. If
\begin{equation}\label{fast}
\limsup\limits_{n\rightarrow\infty}\frac{a_n}{c^n}=\infty \text{  for every positive $c\in {\ensuremath {\mathbb R  }}$,}
\end{equation}
then $\Lambda\subset\Lambda'$ implies that $\Lambda'$ is type 2.
\end{corollary}

This motivates the following definition.
\begin{definition}\label{*supertwo}
{
We say that $\LLL$ is super type $2$ if $\Lambda\subset\Lambda'$  implies that $\Lambda'$ is type $2.$
}
\end{definition}

\begin{question}\label{*qsupertwo}
{\it Find a characterization of super type $2$ sets. Is the typical $\LLL$ super type $2$?}
\end{question}
If this question is too difficult find more necessary
and more sufficient conditions for $\LLL$ being super type $2$.

The next theorem is about sets which are very much 
distant from being super type $2$.
This theorem also demonstrates that the heuristic idea of adding more and more elements to obtain type $2$ sets can fail in a spectacular way.
\begin{theorem}[\cite{[BHMVa]}]{\label{sandwich}}
There exists a collection of discrete sets $\{\Lambda_n\}_{n \in  {\mathbb {Z}}}$ such that
 \begin{equation}\label{sandwichineq}
\Lambda_{n+1} \subset \Lambda_{n} \text{    for all    } n \in  {\mathbb {Z}}
\end{equation}
and $\Lambda_n$ is type 1 if $n$ is odd and type 2 if $n$ is even.
\end{theorem}
 
 Following Theorem 5 in \cite{[BKM1]} it was remarked that
if $\Lambda$ is asymptotically dense and consists of elements
independent over $ {\ensuremath {\mathbb Q}}$ then using Theorem \ref{*BKMth5} it is easy to show that
$\Lambda$ is type $2.$ Recall also Theorem \ref{alg indep type 1} saying that on the other hand there are type $1$ sets containing infinitely many independent elements.

 In the next theorem we see that working with numbers 
 independent over the rationals one can significantly decrease 
 the rate appearing in \eqref{fast} of Corollary \ref{suru 2-es}.
\begin{theorem}[\cite{[BHMVa]}]\label{slow always 2}
If $\Lambda$ consists  of numbers independent over the rationals and
\begin{equation}\label{speed}
\limsup_{n\rightarrow\infty} \frac{\#(\Lambda\cap[0,n))}{n} = \infty
\end{equation}
then every set containing $\Lambda$ is type 2.
\end{theorem}

This can also lead to some questions.
\begin{question}\label{*qrateindep}
{\it What is the best rate in \eqref{speed} with which Theorem \ref{slow always 2}
is still true? For example does ${\log \log n \cdot \#(\Lambda\cap[0,n))}/{n}$, or ${(\log n)^{0.9} \cdot\#(\Lambda\cap[0,n))}/{n}$ work in Theorem \ref{slow always 2}?} 
\end{question}
 
 One may wonder if it is always possible to construct
a chain appearing in Therem \ref{sandwich}
such that $\Lambda_{0}$ is an arbitrary type $2$ set.
Combining Theorems \ref{*BKMth5} and  \ref{slow always 2} we obtain a negative answer:

\begin{theorem}[\cite{[BHMVa]}]\label{always type 2} Assume that $\Lambda$
satisfies
\eqref{speed}
 and $\Lambda$ consists of numbers independent over the rationals.
 In this case for any $\Lambda'$ satisfying \eqref{*proplll} and $\Lambda'\subseteq\Lambda$ or $\Lambda\subseteq\Lambda'$ we have that $\Lambda'$ is type 2. \end{theorem}

Usually adding elements to type $1$ sets can turn them into type $2$ sets, but at least the following is true:
 \begin{proposition}[\cite{[BHMVa]}]\label{*proptos}
If $\Lambda_1,\Lambda_2\subset {\ensuremath {\mathbb R}}$ are type 1 sets then $\Lambda_1\cup\Lambda_2$ is also type 1.
\end{proposition}
 
 The next theorem is a bit of a surprise, but it is in line with Theorem \ref{sandwich}.
 \begin{proposition}[\cite{[BHMVa]}]\label{2-esek unioja}
There exist two decreasing gap asymptotically dense type 2 sets $\Lambda_1$ and $\Lambda_2$ such that their union is type 1.
\end{proposition}

Minkowski sums of sets are much larger than the original ones. Despite of this we still have:
\begin{theorem}[\cite{[BHMVa]}]\label{*minksum}
If the sets $\Lambda=\{\lambda_0,\lambda_1,\ldots\}$ and $\Lambda'=\{\lambda'_0,\lambda'_1,\ldots\}$ are type 1 then the Minkowski sum $\Lambda+\Lambda'$ is
also type 1.
\end{theorem}


The next theorem is the Minkowski sum version of Corollary \ref{suru 2-es}:
\begin{proposition}[\cite{[BHMVa]}]\label{*mink2}
There is a type 2 set $\Lambda$ such that for every $\Lambda'=\{\lambda'_0,\lambda'_1,\ldots\}$ the Minkowski sum $\Lambda+\Lambda'$ is type 2.
\end{proposition}

Analogous questions to the ones asked after Theorem \ref{suru 2-es} can also be stated.
\begin{question}\label{*qsupertwoo}
{\it
Find a characterization (or more sufficient, or more necessary conditions) of sets satisfying the conclusion of Proposition \ref{*mink2}. Does the typical $\LLL$ satisfy this conclusion?}
\end{question}

\section{$\sum_{n=1}^{\infty}c_n f(x+\lambda_n)$}\label{*seccnxpln}

As we mentioned at the end of Section \ref{*secfnxp}
in the periodic case, corresponding to the Khinchin conjecture,
several papers
considered weighted averages like \eqref{*BWa}.
This motivates the following definition:
\begin{definition}[\cite{[BHMVr]}]\label{*defctype}
We say that an asymptotically dense set $\Lambda$ is $\pmb{c}$-type $2$ with respect to the positive sequence $\pmb{c}=(c_n)_{n=1}^{\infty}$,
 if there exists a nonnegative measurable ``witness" function $f$ such that the series $s_{\pmb{c}}(x)=s_{\pmb{c},f}(x)=\sum_{n=1}^{\infty}c_n f(x+\lambda_n)$ does not converge almost everywhere and does not diverge almost everywhere either.
Of course, those $ {\Lambda}$s which are not $\pmb{c}$-type $2$ will be called
$\pmb{c}$-type $1$.
\end{definition}

In the sense of our earlier definition, $\Lambda$ is type $2$ if it is $\pmb{c}$-type $2$ with respect to $c_n\equiv 1$, that is $ {\Lambda}$ is $\pmb{1}$-type $2$. The corresponding convergence and divergence sets
are denoted by $C_{\pmb{c}}(f,\Lambda)$ and $D_{\pmb{c}}(f,\Lambda)$.

The following theorem is a  nice consequence of Theorem \ref{*thwit}.
\begin{theorem}[\cite{[BHMVr]}]\label{*thpmbc1}
 If a set $\Lambda$ is $\pmb{1}$-type $2$, then it is $\pmb{c}$-type $2$ with respect to any positive sequence $\pmb{c}=(c_n)_{n=1}^{\infty}$. \end{theorem}
 
This theorem naturally leads to the next problem.
\begin{question}\label{*vert1}
{\it Characterize (or give more sufficient, or necessary conditions) those $\mathbf{c}$s for which every $\Lambda$ is $\mathbf{c}$-type $2$.}
\end{question}
 
The key property behind Theorem \ref{*thpmbc1} is the fact that for  $\pmb{1}$-type $2$
sets there is always a witness function which is a characteristic function
according to the result of Theorem \ref{*thwit}. This motivates the following definition:
 \begin{definition}\label{*defchi}
A positive sequence $\pmb{c}$ is a $\chi$-sequence
if for any $\pmb{c}$-type  $2$  set $ {\Lambda}$
there is always a characteristic function to witness this property.
\end{definition}

\begin{question}\label{*qchi}
{\it Does Theorem  \ref{*thpmbc1} hold for
all $\chi$-sequences?}
\end{question}

There exist sequences which are not $\chi$-sequences. Indeed, if $\sum c_n$ converges, then for any function $f$ bounded by $K$ we have $s_{\pmb{c},f}(x)=\sum_{n=1}^{\infty}c_{n}f(x+\lambda_{n})\leq \sum_{n=1}^{\infty}c_{n}K$,
and hence $s_{\pmb{c},f}$ converges everywhere. 
On the other hand, by Theorem \ref{*thpmbc1} there are  
$\pmb{c}$-type $2$ sets $\Lambda$, in this case with unbounded
witness functions.

 Hence it is natural to consider the following problem.
\begin{question}\label{*qchica}
{\it 
Characterize $\chi$-sequences.}
\end{question}

The next theorem  shows that not all sequences have this property by showing the other extreme: sequences for which every $ {\Lambda}$ is $\pmb{c}$-type  $2$.

\begin{theorem}[\cite{[BHMVr]}]\label{*univt2seq}
Suppose that $\pmb{c}=(c_n)$ is a sequence of positive numbers satisfying the following condition:
 \begin{equation}\label{cn fast decr}
\sum_{j=n+1}^\infty c_j < 2^{-n}c_n \text{    for every    } n \in  {\mathbb {N}}.
\end{equation}
Then every discrete set $\Lambda$ is $\pmb{c}$-type  $2$.
\end{theorem}


\medskip
\noindent {\bf Acknowledgment}. The author wishes to thank  B. Hanson, B. Maga, and G. Vértesy for their comments on the first draft.


\bibliographystyle{amsplain} 
\bibliography{fnxsurvey}

\end{document}